\documentclass[11pt]{amsart}
\usepackage{amsmath,amssymb,amsfonts,url}
\numberwithin{equation}{section}

\newtheorem{thm}{Theorem}[section]

\begin{document}
\title[Blowup solutions]{Blowup solutions of elliptic systems in two dimensional spaces} 

\author{Lei Zhang}
\address{Department of Mathematics\\
        University of Florida\\
        358 Little Hall P.O.Box 118105\\
        Gainesville FL 32611-8105}
\email{leizhang@ufl.edu}

\maketitle
\date{\today}

 \smallskip

%\setcounter{section}{0}
%%%%%%%%%%%%%%%%%%%%%%%%%%%%%%%%%%%%%%%

Systems of elliptic equations defined in two dimensional spaces with exponential nonlinearity are very commonly observed in geometry, physics and other disciplines of sciences. In this report we survey some results on the following system of equations:
\begin{equation}\label{eq1}
\Delta u_i+\sum_{j=1}^n a_{ij} h_je^{u_j}=\sum_{k=1}^N4\pi \gamma_{ik}\delta_{p_k},\quad B_1\subset \mathbb R^2,\quad i=1,...,n
\end{equation}
where $A=(a_{ij})_{n\times n}$ is the coefficient matrix, $h_i$ are positive smooth functions and, $\gamma_{ik}>-1 $ are constants, $\delta_{p_k}$ represents Dirac mass at $p_k$, $B_1$ is the unit ball in $\mathbb R^2$. There are different variations of (\ref{eq1}) depending on the context and applications. If the system is defined on a Riemann surface $M$, (\ref{eq1}) is often written in the following form:
$$\Delta_g u_i +\sum_{j=1}^n a_{ij} h_j e^{u_j} -K(x)= 4\pi \sum_{k=1}^N \gamma_{ik} \delta_{p_k} $$
where $K$ is the Gauss curvature.

If $A$ is nonnegative (i.e. $a_{ij}\ge 0$ for all $i,j$), (\ref{eq1}) is called a Liouville system. If some entries of $A$ are negative, the most famous case is the Toda system, where $A$ is the following Cartan matrix:
$$A=\left (\begin{array}{ccccc}
2 & -1 & 0 & ... & ... \\
-1& 2 & -1 & ... & ... \\
0 & -1 & 2 & ... & ... \\
... & ... & ... & ... & ... \\
... & ... & ... & -1 & 2
\end{array}
\right )
$$
If the system (\ref{eq1}) has only one equation, then it is reduced to the single Liouville equation
$$\Delta u+ h e^u=0,\quad B_1.  $$
The main purpose of this report is to compare the behavior of blowup solutions to the Liouville equation, Liouville systems and Toda systems. The comparison will be focused on three perspectives: The classification of global solutions, concentration of energy and application of the corresponding equation/system defined on Riemann surfaces.

\bigskip

In 1992 Chen and Li \cite{chen-li-1} proved that all solutions to
\begin{equation}\label{liou1}
\Delta u+ e^u=0,\quad \mbox{ in } \mathbb R^2,\quad \int_{\mathbb R^2}e^u<\infty
\end{equation}
can be explicitly written as
$$u_{\lambda, x_0}(x)=\log
\frac{\lambda^2}{(1+\frac{\lambda^2}8 |x-x_0|^2)^2},\quad \lambda>0,\quad
x_0\in \mathbb R^2 $$
where $\lambda>0$ and $x_0\in \mathbb R^2$ are three parameters that describe all the solutions as a family. It is easy to verify that all solutions in this family have one energy: $\int_{\mathbb R^2} e^u=8\pi$. If the right hand side of (\ref{liou1}) has a Dirac mass at $0$: $4\pi \gamma \delta_0$. Then if $\gamma>-1$, Prajapat and Tarantello \cite{prajapat} also classified all the global solutions. From the explicit expressions of global solutions to the singular Liouville equation we can see that there is also one energy $\int_{\mathbb R^2} e^u=8\pi(1+\gamma)$ for all solutions, a fact first discovered by Chen-Li \cite{chen-li-duke2} before the classification theorem.

For Liouville systems defined on $\mathbb R^2$:
$$\Delta u_i+\sum_{j=1}^n a_{ij} e^{u_j}=0,\quad \mathbb R^2, \quad i=1,...,n, \quad \int_{\mathbb R^2} e^{u_i}<\infty, $$
the following assumption is more or less standard:
$$(H_1): \quad A \mbox{  is symmetry, invertible and irreducible}. $$
The coefficient matrix $A$ being irreducible means the system can not be written as two separate smallers systems. In 1997 Chipot-Shafrir-Wolansky \cite{chipot-1} established a mapping from the set of global solutions $u=(u_1,...,u_n)$ to the set of energies: $\sigma=(\sigma_1,...,\sigma_n)$ (where $\sigma_i=\frac{1}{2\pi}\int_{\mathbb R^2} e^{u_i}$). They proved, among other things, that
\begin{thm}\label{chipot-t1}(Chipot,Shafrir,Wolanski 97) All components of $u$ are radial with respect to a common point. Also ( for $I=\{1,...,n\}$)
$$\Lambda_I(\sigma)=4\sum_{i=1}^n\sigma_i-\sum_{i,j=1}^n a_{ij}\sigma_i\sigma_j=0,\quad \Lambda_J>0, \forall \emptyset\varsubsetneqq J\varsubsetneqq I $$
On the other hand, for any $\sigma$ satisfying $\Lambda_I=0$ and $\Lambda_J>0$ for all $\emptyset\varsubsetneqq J\varsubsetneqq I$, there is a global solution corresponding to it.
\end{thm}
One major difference between Liouville systems and Liouville equation is on the energy. From Theorem \ref{chipot-t1} we see that the energy of Liouville system is a $n-1$ dimensional hypersurface. So there is a continuum of energy for Liouville systems. Theorem \ref{chipot-t1}
is not the classification result for global solutions of Liouville system since it does not address the uniqueness property of the mapping from global solutions to the hypersurface of energy. The uniqueness property is solved by Lin and Zhang in 2010:
\begin{thm}\label{lin-zhang-t1} (Lin-Zhang, 10)
Let $u=(u_1,..,u_n)$ and $v=(v_1,..,v_n)$ be entire solutions such that $\int_{\mathbb R^2} e^{u_i}=\int_{\mathbb R^2} e^{v_i}$. Then
$$v_i(y)=u_i(\delta y+x_0)+2\log \delta $$ for some $x_0\in \mathbb R^2$ and $\delta>0$.
\end{thm}
The combination of Theorem \ref{chipot-t1} and Theorem \ref{lin-zhang-t1} gives the classification of all global solutions of Liouville systems under $(H_1)$.

As far as the classification of the global solutions of Toda system is concerned, if the system has one possible singularity:
 $$\Delta u_i+\sum_{j=1}^n a_{ij} e^{u_j} =4\pi \gamma_i \delta_0, \,\, \mbox{ in } \mathbb R^2,
 \int_{\mathbb R^2}e^{u_i}<\infty \quad i=1,...,n,$$
 where $A=(a_{ij})_{n\times n}$ is the Cartan matrix, $\gamma_i>-1$.
 A recent result of Lin-Wei-Ye \cite{lwy} completely answered this question and they proved that: 1. all the solutions to $SU(n+1)$ Toda systems with one singular source can be written explicitly as a family of $n^2+2n$ parameters, 2. For fixed $\gamma_i$, all global solutions have one energy: suppose $u=(u_1,...,u_n)$ is a global solution, then
$$\sum_{j=1}^n a_{ij} \int_{\mathbb R^2} e^{u_j} dx = 4\pi (2+\gamma_i+\gamma_{n+1-i}),\quad i=1,...,n. $$
The classification theorem of global solutions to $SU(n+1)$ Toda system without singular source was proved by Jost-Wang in \cite{jost-wang}.

\medskip

The classification of global solutions has great influence on the concentration of energy of blowup solutions near its blowup points. For Liouville equations we take the following standard case as an example: Let $u_k$ be a sequence of blowup solutions to
$$\Delta u_k + h e^{u_k}=0, \quad B_1\subset \mathbb R^2 $$
where $h$ is a positive smooth function. Suppose there is an upper bound on the energy: $\int_{B_1} h e^{u_k}<C$. Brezis-Merle \cite{bm} proved that when blowup phenomenon occurs, $h e^{u_k}\rightharpoonup \sum_l \alpha_l \delta_{p_l}$ for some $p_1,...,p_N\in B_1$ and $\alpha_l\ge 4\pi$. Then Li-Shafrir \cite{li-shafrir} proved that each $\alpha_i$ is a multiple of $8\pi$. The reason is around each blowup point $p_l$, there are a few sequences of small disks all tending to the blowup point. These small disks are all shrinking to the blowup point and are mutually disjoint. The profile of $u_k$ in each small disk is roughly like that of the global solution. In other words the energy of $u_k$ in each small disk is roughly equal to $8\pi$ and there is little energy outside these disks. For a while people wondered if there is only one sequence of disks tending to each blowup point. This is not the case, X. Chen \cite{xxchen} proved that even for $h\equiv 1$ any multiple of $8\pi$ may occur as the weak limit of the energy. In 2000 Li \cite{licmp} proved that if the blowup solutions are assumed to have \emph{bounded oscillation on $\partial B_1$} (suppose $0$ is the only blowup point in $B_1$) then the blowup picture is much simpler: $h e^{u_k}\rightharpoonup 8\pi \delta_0$ and
$$|u_k(x)-\log \frac{e^{u_k(0)}}{(1+\frac{h(0)}8 e^{u_k(0)}|x|^2)^2}|=O(1) \mbox{ in } B_1. $$
Note that the second term on the left hand side is the global solutions with appropriate scaling. We shall address Li's estimate as $O(1)$ estimate.
The \emph{bounded oscillation assumption } around a blowup point is a natural assumption since in most applications it is automatically satisfied.

Li's $O(1)$ estimate is very useful. In particular it implies that blowup solutions look roughly the same around different blowup points and the energy outside the bubbling area is very small. The $O(1)$ estimate for singular Liouville equations was proved by Bartolucci-Chen-Lin-Tarantello \cite{bclt}.

The $O(1)$ estimate for systems is much harder since for systems an essential difficulty is the ``partial blown-up phenomenon", which means when the blowup solutions are scaled according to the maximum of all components, some components may disappear in the limit. In \cite{linzhang2} Lin and Zhang proposed a second assumption on the coefficient matrix $A$ of Liouville systems ( in addition to $(H_1)$):
$$(H_2):\quad a^{ii}\leq 0 \,\, \forall i, \quad a^{ij}\geq 0,\,\, \forall i\neq j,\quad
\sum_{j=1}^na^{ij}\geq 0, \,\, \forall i $$
where $a^{ij}$s are the entries in $A^{-1}$.

The assumption $(H_2)$ is a strong interaction assumption. If $A=(a_{ij})_{n\times n}$ is $2\times 2$ matrix, the combination of $(H_1)$ and $(H_2)$ is equivalent to
$$a_{ij}\ge 0,\quad \max\{a_{11},a_{22}\}\le a_{12},\quad det(A)\neq 0. $$
Since $a_{12}$ is the connection between two equations, $(H_2)$ clearly indicates that the ties between two equations must be strong.  Lin-Zhang
\cite{linzhang2} proved that under $(H_1)$ and $(H_2)$ the profile of blowup solutions around different blowup points must be $O(1)$ close, even when partial blown-up phenomenon occurs.

The $O(1)$ estimate for Toda system is much harder. Much less has been obtained so far even for $SU(3)$ Toda system ( the system that has only two equations). In \cite{jost-lin-wang} Jost-Lin-Wang proved the $O(1)$ estimate for $SU(3)$ Toda systems under one essential assumption: The blowup sequence has to converge to a full system ( of two equations) after scaling. Since the Toda system is an integrable system and is so much involved with algebraic geometry, Jost-Lin-Wang's proof actually uses tools from holonomy theory. Recently Lin-Wei-Zhang \cite{lwz-2} proved a better than $O(1)$ estimate for fully blown-up sequence of general $SU(n+1)$ Toda systems (no restriction on the number of equations) without using tools from algebraic geometry. Without assuming the fully blown-up picture the behavior of blowup solutions is very complicated. Jost-Lin-Wang proved that for $SU(3)$ Toda systems, the energy concentration can only be one of the following six types near a blowup point:
$$(0,0),(0,4\pi),(4\pi,0),(4\pi, 8\pi),(8\pi,4\pi),(8\pi,8\pi). $$
It is not clear whether each of the six types only corresponds to one asymptotic behavior of solutions. However Pistoia-Musso-Wei \cite{pmw} proved that all six types are possible. Recently Lin-Wei-Zhang \cite{lwz-1} extended the result of Jost-Lin-Wei to singular $SU(3)$ Toda systems.

Even though the $O(1)$ estimate provides precise asymptotic behavior of blowup solutions near a blowup point, more accurate estimates are needed for some applications. For Liouville equation of the mean field type
\begin{equation}\label{mf}
\Delta_g u + \rho (\frac{he^u}{\int_Mhe^udV_g}-1)=0
\end{equation}
we assume that the volume of $M$ is $1$ and the average of $u$ is $0$. Then the topological degree of $u$ is well defined if $\rho<8\pi$ or  $8N\pi
<\rho<8(N+1)\pi$ ($N\in \mathbb N$). Chen-Lin \cite{chen-lin-sharp,chenlin2} proved the following degree counting formula:
\begin{equation}\label{deg1}
d_{\rho}=\left\{\begin{array}{ll}
1 \quad \rho<8\pi,\\ \\
\frac{(-\chi_M +1)...(-\chi_M +N)}{N!}\quad 8N\pi<\rho<8(N+1)\pi.
\end{array}
\right.
\end{equation}
and their key estimate is much stronger than $O(1)$ estimate. From (\ref{deg1}) we observe that if $\chi_M\le 0$, that is, the genus of $M$ is great than $0$, the mean field equation (\ref{mf}) always has a solution if $\rho$ is not a multiple of $8\pi$. Another degree counting formula for the Liouville equation with singular source was also established by Chen-Lin \cite{chen-lin-sing}.

It is easy to imagine that the degree counting formula for systems is generally much harder, since the blowup picture is lot more complicated. For Liouville systems defined on Riemann surfaces:
$$
\Delta_g u_i+\sum_{j=1}^n\rho_j a_{ij} (\frac{h_je^{u_j}}{\int_Mh_je^{u_j}dV_g}-1)=0,\quad i=1,..,n
$$
where $h_1,...,h_n$ are positive smooth functions, $\rho_1,..,\rho_n\ge 0$ are constants, the average of each $u_i$ is $0$. If the coefficient matrix satisfies $(H_1)$ and $(H_2)$, Lin and Zhang \cite{linzhang2} proved that $u=(u_1,...,u_n)$ satisfies the a priori estimate $|u_i|\le C$ if $\rho=(\rho_1,...,\rho_n)$ is not on the following family of critical hyper-surfaces:
\begin{equation*}
\Gamma_N=
\left\{
\rho\,\big|\, 8\pi N \sum_{i=1}^n\rho_i=\sum_{i,j}a_{ij}\rho_i\rho_j
\right\}.
\end{equation*}
The Leray-Schauder degree of the solution is well defined on regions bounded by these hyper-surfaces in the first hyper-octant. Let the region below $\Gamma_1$ be called $\mathcal{O}_0$, the region between $\Gamma_N$ and $\Gamma_{N+1}$ be called $\mathcal{O}_N$. Lin-Zhang 
\cite{linzhang2} proved the following degree counting formula:
\begin{equation*}
d_{\rho}=\left\{\begin{array}{ll}
1\quad \text{if }\rho \in \mathcal{O}_0\\
\frac{1}{N!}\bigg ((-\chi_M+1)...(-\chi_M+N) \bigg )\quad \text{if }\rho \in \mathcal{O}_N.
\end{array}
\right.
\end{equation*}
where $\chi_M$ is the Euler characteristic of $M$. Since $\chi_M=2-2g_e$ ($g_e$ is the genus of $M$), the degree is always positive if the genus is greater than $0$ and $\rho\not \in \Gamma_N$.

For Toda systems a degree counting theorem is not yet established. One main reason is the blowup picture is very complicated, as one can see from the work of Jost-Lin-Wang \cite{jost-lin-wang} and Lin-Wei-Zhang \cite{lwz-1}. In \cite{lwz-1} it is proved that for $SU(3)$ Singular Toda systems there are only finite ways for the energy of blowup solutions to concentrate. This result leads to the recent existence result of Battaglia-Jevnikar-Malchiordi-Ruiz \cite{andrea1}. However for Liouville systems, since the energy belongs to a $n-1$ dimensional hyper-surface, major difficulties appear when a very precise comparison of bubbling solutions around different blowup points is needed. Regarding Lin-Zhang's degree counting theorem for Liouville systems, if $\rho$ is in a compact subset of $\mathcal{O}_N$, there is a priori estimate on the solution $u$. If $\rho$ tends to $\Gamma_N$, the blowup phenomenon occurs and there are exactly $N$ blowup points. The abundance of the energy for Liouville systems causes essential difficulty to the study of the interaction of bubbling solutions. When $N=1$, there is only one blowup point when $\rho\to \Gamma_1$ ( so there is no need to study bubble interaction in this simple case), the study of this case already reveals very surprising information, as Lin-Zhang \cite{linzhang3} found out that there is one special point $q$ on $\Gamma_1$ such that if $\rho\to q$ the leading term of $\rho-q$ only involves the curvature information at the corresponding blowup point. However if $\rho$ tends to any other point on $\Gamma_1$, the distance between $\rho$ and $\Gamma_1$ is an integral which is involved with the geometry of the whole manifold. The reader may look into \cite{linzhang3} for the precise statement.

\end{document}